\documentstyle{article}

\def\QED{\quad\blackslug\lower 8.5pt\null}

\begin{document}

\title{
{\large \bf A CONFORMAL DIFFERENTIAL INVARIANT 
  AND THE CONFORMAL RIGIDITY 
 OF HYPERSURFACES}
}
\renewcommand{\footnote}{}
\author{
 M.A. Akivis \and  V.V. Goldberg
\footnotetext{{\em $1991$ 
Mathematics Subject Classification}:  Primary 53A30. \newline  
\hspace*{5mm}{\em Keywords and phrases}: 
Conformal and pseudoconformal geometry, 
hypersurface, first and second fundamental forms, conformal 
quadratic element, moving frames, conformal rigidity. \newline 
\hspace*{5mm}This research was   partially supported by Volkswagen-Stiftung 
(RiP-program at MFO). The research of the first author was 
also partially supported by the Israel 
Ministry of Absorption and the Israel Public Council for Soviet 
Jewry. 
} 
}
\date{}

\maketitle


{\begin{abstract}  
\noindent
For a hypersurface $V^{n-1}$  
of a conformal space,  we introduce a conformal differential 
invariant $I = \frac{h^2}{g}$, where $g$ and $h$ are 
the first  and the second 
fundamental forms of  $V^{n-1}$ connected by the apolarity condition. 
This invariant is called the {\em conformal quadratic element} 
of $V^{n-1}$. The solution of the problem of conformal 
rigidity is presented  in the framework  of conformal 
differential  geometry and connected with the conformal 
quadratic element of $V^{n-1}$.   
The main theorem states: 

  Let $n \geq 4$ and $V^{n-1}$ and $\overline{V}^{n-1}$ be two nonisotropic 
hypersurfaces without umbilical points 
in a conformal space $C^n$ or a 
pseudoconformal space $C^n_q$ of signature $(p, q), \;\; p = 
n - q$. Suppose that there is a one-to-one correspondence 
$f: V^{n-1} \rightarrow \overline{V}^{n-1}$ between points of
 these hypersurfaces, and in the corresponding points of 
$V^{n-1}$ and $\overline{V}^{n-1}$ the following condition holds:
$
\overline{I} = f_* I, 
$
 where $f_*: T (V^{n-1}) \rightarrow T (\overline{V}^{n-1})$ 
is a mapping induced by the correspondence $f$. 
 Then   the hypersurfaces $V^{n-1}$ and 
$\overline{V}^{n-1}$ are conformally equivalent.
\end{abstract}
\setcounter{equation}{0}

\vspace*{5mm}
{\bf 1.} 
 In local differential geometry the rigidity theorems contain 
conditions under which two submanifolds of a homogeneous space 
can differ only by a location in the space. For hypersurfaces 
in a projective space, the rigidity problem was considered 
by G. Fubini [F 16, 18] (see also pp. 605--629 of
 the book [F\v{C} 26] by Fubini and E. \v{C}ech), 
\'{E}. Cartan [C 20], G. Jensen and E. Musso [JM 94], 
and by the authors of this paper in the book [AG 93] 
(Section {\bf 7.4}). 

The problem of conformal rigidity of submanifolds is also of 
great interest. This problem was studied  by Cartan [C 17],  
M. do Carmo and M. Dajczer [CD 87] and R. Sacksteder [S 62] 
(see also the paper [Su 82] by R. Sulanke in which the author 
considered problems close to the rigidity problem). 
However, in these papers the problem of conformal rigidity was 
investigated in the framework of Euclidean geometry. 

In the current paper we present the solution of this problem 
in the framework of conformal differential  geometry. 
To this end, we introduce a conformal quadratic element and prove that 
if $n \geq 4$ and there exists a one-to-one point 
correspondence of two hypersurfaces both not having umbilical 
points preserving this quadratic element, then 
the hypersurfaces are conformally equivalent. 
Moreover, 
we consider the  rigidity  problem not only for hypersurfaces of 
a conformal space but also  for hypersurfaces of a 
pseudoconformal space. We only assume that a hypersurface is not 
isotropic, i.e.  its tangent subspaces are not tangent to the 
isotropic cones.

{\bf 2.} Let $V^{n-1}$ be a nonisotropic  hypersurface of a 
real conformal space $C^n$ or a real pseudoconformal space
 $C^n_q$ of signature $(p, q)$, where $p = n - q$. With any point 
$x \in V^{n-1}$, we associate a conformal moving frame consisting 
of two points $A_0 =x$ and $A_{n+1}$ and $n$ independent 
hyperspheres $A_1, \ldots , A_n$, passing through these two 
points. We will assume that the hypersphere $A_n$ is tangent to 
the  hypersurface $V^{n-1}$ at the point $A_0$, and the 
hyperspheres $A_i, i = 1, \ldots, n - 1$, are orthogonal to 
$V^{n-1}$  at the point $A_0$. Then the frame elements satisfy 
the following conditions:
\begin{equation}\label{eq:1}
\renewcommand{\arraystretch}{1.3}
\begin{array}{ll}
 (A_0, A_0) = (A_{n+1}, A_{n+1}) = 0, &  (A_0, A_i) =  
(A_{n+1}, A_i) = 0, \\
 (A_0, A_n) = (A_{n+1}, A_n) = 0, &  (A_i, A_n) = 0,
\end{array}
\renewcommand{\arraystretch}{1}
\end{equation}
where parentheses denote the scalar product of corresponding 
frame elements. The first two of these conditions mean that the 
frame elements $A_0$ and $A_{n+1}$ are points, the following four 
conditions mean that  the hyperspheres $A_i$ and $A_n$ pass 
through these two points, and finally, the last condition 
expresses the orthogonality of  the hyperspheres $A_i$ and $A_n$. 
In addition, we normalize the points  $A_0$ and $A_{n+1}$ by the 
condition 
\begin{equation}\label{eq:2}
 (A_0,  A_{n+1}) = -1.
\end{equation}
We will not demand the orthogonality of the hyperspheres $A_i$ 
and will write their scalar products in the form:
\begin{equation}\label{eq:3}
 (A_i,  A_j) = g_{ij}.
\end{equation}
where $\det (g_{ij}) \neq 0$, since the hypersurface $V^{n-1}$ is 
not isotropic.

If $X$ is an arbitrary point of  the space $C^n$ or the space  
$C^n_q$, then it can be represented as 
$$
X = x^0 A_0 + x^i A_i + x^n A_n + x^{n+1} A_{n+1}. 
$$
Since for any point $X \in C^n$, we have $(X, X) = 0$, then it 
follows from this that the coordinates of the point $X$ satisfy 
the equation 
$$
g_{ij} x^i x^j + \epsilon (x^n)^2 - 2 x^0 x^{n+1} = 0, 
$$
where $\epsilon = (A_n, A_n) \neq 0$, since the quadratic 
form on the left-hand side 
of the last equation is nondegenerate.
The last equation  is the equation of a nondegenerate hyperquadric 
$Q^n$ of a projective space $P^{n+1}$ onto which the space $C^n$  
or the space  $C^n_q$ are mapped under the Darboux mapping. The 
left-hand side of this equation is of signature $(p+1, q+1)$. 
Under the Darboux mapping, the images of the points $A_0$ and 
$A_{n+1}$ are the points lying on the hyperquadric $Q^n$, and 
in general, the 
images of the hyperspheres $A_i$ and $A_n$ are the points not 
belonging to $Q^n$.

We will prove now that the quantity $\epsilon$ can be always reduced 
to 1. In fact, if $\epsilon > 0$, then this can be achieved by 
means of renormalization of the hypersphere $A_n$. If $\epsilon <0$, 
then we can replace the point $A_0$ by $- A_0$ which implies 
$$
g_{ij} x^i x^j + \epsilon (x^n)^2 + 2 x^0 x^{n+1} = 0,
$$
and then change the sign of the left-hand side of the above 
equation and again reduce to 1 the positive quantity $-\epsilon$. 
As a result, the equation of the Darboux hyperquadric takes the form 
$$
- g_{ij} x^i x^j +  (x^n)^2 - 2 x^0 x^{n+1} = 0.
$$
By setting $g_{ij} = - \widetilde{g}_{ij}$ and suppressing tilde, 
we reduce the last equation to the form
\begin{equation}\label{eq:4}
 g_{ij} x^i x^j +  (x^n)^2 - 2 x^0 x^{n+1} = 0.
\end{equation}

Therefore, for any $\epsilon \neq 0$, we can normalize the hypersphere 
$A_n$ in such a way that 
\begin{equation}\label{eq:5}
(A_n, A_n) = 1.
\end{equation}
The form $g_{ij} x^i x^j$ in equation (4) has signature $(p - 1, q)$.

{\bf 3.}  The equations of infinitesimal displacement of our 
conformal frame have the form:
\begin{equation}\label{eq:6}
dA_\xi = \omega_\xi^\eta A_\eta, \;\; \xi, \eta = 0, 1, \ldots ,  n + 1,
\end{equation}
where $\omega_\xi^\eta$ are 1-forms containing parameters, on  
which the family of frames in question depends, and their 
differentials: $\omega_\xi^\eta = \omega_\xi^\eta (u, du)$. 
By equations (1)---(3) and (5), these forms must satisfy the following 
conditions:
\begin{equation}\label{eq:7}
\renewcommand{\arraystretch}{1.3}
\left\{
\begin{array}{ll}
\omega_0^{n+1} = \omega_{n+1}^0 = 0, &
 \omega_0^0 + \omega^{n+1}_{n+1} = 0,\\
\omega_i^{n+1} = g_{ij} \omega_0^j, &  
 \omega_{n+1}^i = g^{ij} \omega_j^0,\\ 
\omega_n^{n+1} = \omega_0^n, & 
\omega^n_{n+1} = \omega_n^0, \\
 \omega_i^n = - g_{ij} \omega_n^j, & 
\omega_n^n = 0,\\
dg_{ij} = g_{ik} \omega_j^k + g_{kj} \omega_i^k. & 
\end{array}
\right.
\renewcommand{\arraystretch}{1}
\end{equation}
In addition, the forms $\omega_\xi^\eta$ satisfy the structure 
equations of the conformal space:
\begin{equation}\label{eq:8}
d \omega_\xi^\eta = \omega_\xi^\zeta \wedge \omega_\zeta^\eta, 
\;\;\; \xi, \eta, \zeta = 0, 1, \ldots , n + 1,
\end{equation}

{\bf 4.} Since the hypersphere $A_n$ is tangent to the 
hypersurface $V^{n-1}$  at the point $x = A_0$, the condition
 $(A_n, dA_0) = 0$ holds. It follows from this condition that
\begin{equation}\label{eq:9}
 \omega_0^n = 0.
\end{equation}
This equation determines the family of frames of first order 
associated with the hypersurface $V^{n-1}$. This family can be 
considered as the frame bundle ${\cal R}^1 (V^{n-1})$  of first order 
 with the base $V^{n-1}$. Its fiber 
is a set of frames  which is associated with the point 
$x \in V^{n-1}$ in the manner indicated above. 
The structure group of the frame bundle ${\cal R}^1 (V^{n-1})$ 
is a subgroup of the fundamental group of 
the space $C^n_q$ whose transformations 
leave invariant the tangent element of 
$V^{n-1}$ consisting of a point $x \in V^{n-1}$ and the 
tangent subspace $T_x (V^{n-1})$. The 1-forms $\omega_0^i$,  
which we will denote further by  $\omega^i$, are basis forms of 
the frame bundle ${\cal R}^1 (V^{n-1})$, and the 1-forms 
$\omega_0^0, \omega_i^0, \omega_n^0$ and $\omega_j^i$ are its 
fiber forms.

By equation (9), on the hypersurface $V^{n-1}$, we have 
\begin{equation}\label{eq:10}
d A_0 =  \omega_0^0 A_0 + \omega^i A_i.
\end{equation}
By virtue of this, 
$$
(d A_0, d A_0) = g_{ij} \omega^i \omega^j.
$$
The quadratic form $g =  g_{ij} \omega^i \omega^j$ is 
relatively invariant and determines a conformal structure on 
the hypersurface $V^{n-1}$. This form is nondegenerate and of 
signature $(p-1, q)$. For $q = 0$, i.e. for the proper conformal 
space $C^n$, the form $g$ is positive definite. For $p = 1$ 
and $q = n - 1$, i.e. for a conformal space $C^n_1$ of Lorentzian 
signature, there exists hypersurfaces with signature $(0, q)$. 
Such hypersurfaces are called {\em spacelike}, and a conformal structure 
induced on them is properly conformal.  The equation 
$ g_{ij} \omega^i \omega^j = 0$ determines the isotropic cone 
of the hypersurface $V^{n-1}$.

Taking the exterior derivative of equation (9) and applying (8), 
we obtain:
\begin{equation}\label{eq:11}
\omega^i \wedge \omega_i^n = 0.
\end{equation}
Applying Cartan's lemma to equation (11), we find that 
\begin{equation}\label{eq:12}
 \omega_i^n = \lambda_{ij} \omega^j, \;\;  \lambda_{ij}  
=  \lambda_{ji}. 
\end{equation}

Differentiating equation (10), we obtain
$$
\renewcommand{\arraystretch}{1.3}
\begin{array}{ll}
d^2 A_0 = &(d \omega_0^0 + (\omega^0_0)^2 
+ \omega^i \omega_i^0) A_0 
+ (\omega_0^0 \omega^i +  \omega^j \omega_j^i) A_i\\ 
&+ \omega^i \omega_i^n A_n  +  \omega^i \omega_i^{n+1} A_{n+1}.
\end{array}
\renewcommand{\arraystretch}{1}
$$ 
It follows that 
$$
\overline{h} = (d^2 A_0, A_n) = \omega^i \omega_i^n = \lambda_{ij} \omega^i \omega^j.
$$
In the tangent subspace $T_x (V^{n - 1})$, 
the equation $\overline{h} = 0$ determines the cone of directions 
along which the hypersphere $A_n$ has a second order tangency 
with $V^{n-1}$.

But the quadratic form $\overline{h}$ is not conformally 
invariant since all hyperspheres of the pencil 
$$
A_n' = A_n + s A_0
$$
as well as the hypersphere $A_n$ are tangent to the hypersurface 
$V^{n-1}$. In view of this, we obtain the following  this pencil 
of  quadratic forms:  
$$
(d^2 A_0, A_n') 
= \omega^i \omega_i^n - s \omega^i \omega_i^{n+1} 
= (\lambda_{ij} - s g_{ij}) \omega^i \omega^j.
$$
and all quadratic forms of this pencil have equal rights. 
From this pencil we will distinguish one form:
$$
h = h_{ij} \omega^i \omega^j,
$$
where
$$
h_{ij} = \lambda_{ij} - \lambda g_{ij} \;\; \mbox{and}\;\; 
\lambda =\frac{1}{n-1} \lambda_{ij} g^{ij}.
$$
It is easy to see that the coefficients $h_{ij}$ satisfy 
the apolarity condition (the trace-free condition):
\begin{equation}\label{eq:13}
 h_{ij} g^{ij} = 0,   
\end{equation}
where $g^{ij}$ is the inverse tensor for $g_{ij}$.

One can prove that geometrically condition (13) means  that 
the cone determined by the equation $h = 0$ is real and there exists 
an orthogonal $(n-1)$-hedron 
formed by tangent directions to $V^{n-1}$ at the point $x$  
 which is inscribed into this cone  (see [AG 93], pp. 214--216).

The tangent hypersphere $C_n = A_n + \lambda A_0$ is conformally invariant. 
It is called the {\em central tangent hypersphere}. 
The cone $h = 0$ is composed of directions 
along which the hypersphere $C_n$ has a second order tangency with 
$V^{n-1}$. 

The quadratic forms $g$ and $h$ are called the {\em first and 
second fundamental forms} of 
the hypersurface $V^{n-1}$.

Points of 
the hypersurface $V^{n-1}$,  
in which the forms $\overline{h}$ and $g$ are proportional, 
are called {\em umbilical points}. It is well-known that if
 at any point $x \in V^{n-1}$, we 
have $\overline{h} = \kappa g$, then 
the hypersurface $V^{n-1}$ 
is a hypersphere or its open part. Note that by (13) the 
condition $\overline{h} = \kappa g$ implies $h = 0$.

If in the frame ${\cal R}_x (V^{n-1})$ associated with 
a point $x \in V^{n-1}$ we replace the tangent hypersphere $A_n$ 
by the central tangent hypersphere $C_n$, then we obtain 
a second order frame ${\cal R}_x^2 (V^{n-1}) \subset {\cal R}^2 
(V^{n-1})$. Transformations of the structural group of the fibre bundle 
${\cal R}^2 (V^{n-1})$ leave invariant not only the tangent element 
of the hypersurface $V^{n-1}$ but also the central tangent 
hypersphere $C_n$ attached to a point $x \in V^{n-1}$. 
In the fiber bundle ${\cal R}^2 (V^{n-1})$, the number 
of fiber forms will be reduced, since the 1-form $\omega_n^0$ 
becomes a linear combination of the basis forms $\omega^i$. 

With respect to a second order frame, equation (12) takes 
the form 
\begin{equation}\label{eq:14}
 \omega_i^n = h_{ij} \omega^j. 
\end{equation}

The point $A_0 \in V^{n-1}$ admits the normalization 
\begin{equation}\label{eq:15}
A_0' = r A_0, \;\; r \neq 0.
\end{equation}
In order to preserve condition (2), we also normalize the point 
$A_{n+1}$ as follows: $A_{n+1}' = \frac{1}{r} A_{n+1}$. 
Since under this normalization we have
$$
(d A_0', d A_0') = r^2 (dA_0, dA_0),
$$
the quadratic form $g$ undergoes the transformation
\begin{equation}\label{eq:16}
g' = r^2 g.
\end{equation}
Moreover, we have
$$
(d^2 A_0', C_n) = r (d^2 A_0, C_n),
$$
and as a result,
\begin{equation}\label{eq:17}
h' = r h.
\end{equation} 
It follows from relations (16) and (17) that {\em the fundamental 
forms $g$ and $h$ of the hypersurface $V^{n-1}$ 
are relatively invariant forms of weight 2 and 1, respectively}.

The forms $g$ and $h$ allow us to construct the expression 
\begin{equation}\label{eq:18}
I = \frac{h^2}{g},
\end{equation}
which is invariant with respect  to normalization (15).  
Hence, this expression is an absolute conformal invariant 
of the hypersurface $V^{n-1}$. It is determined in a  second 
differential neighborhood of $V^{n-1}$.  
The expression $I = I (x, \omega^i)$ is a 
homogeneous  function of second order with respect to the 
basic forms $\omega^i$ of 
the hypersurface $V^{n-1}$. 
We will call this function  
the {\em  conformal quadratic element} of  
the hypersurface $V^{n-1}$ (cf. with the projective 
linear element of $V^{n-1} \subset P^n$ considered in 
[F 16], [F 18], [F \v{C} 26], [JM 94], [AG 93]). 
At   umbilical points, the invariant $I$ vanishes, 
since $h = 0$ at these points.

Since the invariant $I$ is defined by means of the 
first and second quadratic fundamental forms  
the hypersurface $V^{n-1}$, {\em it is a second order 
invariant}. 

\vspace*{2mm} 
\noindent
{\bf Lemma} {\em For $n \geq 4$, at  unumbilical points of 
the hypersurface $V^{n-1}$,  
the   conformal quadratic element $I$ is not 
a quadratic form with respect to $\omega^i$.}

{\sf Proof.} Suppose that 
$$
\frac{h^2}{g} = \theta,
$$
where $\theta = \theta_{ij} \omega^i \omega^j$ is a quadratic form, i.e. 
\begin{equation}\label{eq:19}
h^2 = g \theta.
\end{equation}
 
Suppose that $\mbox{rank}\; h = \rho$. We have $\rho \neq 0$, 
since the points under consideration are not umbilical. 
Moreover, we have $\rho \neq 1$, since if $\rho = 1$, then 
there exists a coordinate system in which all components 
of the tensor $h_{ij}$ vanish except $h_{11}$. But then 
condition (13) reduces to $h_{11} g^{11} = 0$. Since $g^{11} 
\neq 0$, it follows that $h_{11} = 0$ and $\rho = 0$.
Consider the cases $\rho = 2$ and $\rho > 2$ 
 separately.  For $\rho = 2$, 
 the form $h$ can be always written as
$$
h = \alpha \cdot \beta,
$$
where $\alpha$ and $\beta$ are linear forms with respect 
to $\omega^i$, since by (13) the form $h$ is alternating. 
Thus, from equation (19) it follows that 
the quadratic form $g$ is also decomposable. However, this is 
impossible since the form $g$ is nondegenerate and for $n \geq 4$, 
is not decomposable into the product of linear factors.

If $\rho > 2$, then the forms $g$ and $h$ are not decomposable, 
and hence equation (19) is possible only if $g$ and $h$ 
are proportional, i.e. only at  umbilical points. 
\rule{3mm}{3mm}

{\bf 5.} 
Consider two smooth, oriented, connected and simply connected 
 hypersurfaces $V^{n-1}$ 
and $\overline{V}^{n-1}$ of a  
conformal space $C^n$ or a pseudoconformal space $C^n_q$. 
Suppose that there is a one-to-one 
correspondence $f: V^{n-1} \rightarrow \overline{V}^{n-1}$ 
under which $f (x) = \overline{x}$, where $x \in V^{n-1}$ and 
$\overline{x} \in \overline{V}^{n-1}$. 
 The correspondence $f$ induces a 
mapping $f_*$ of the tangent bundle $T (V^{n-1})$  
onto the tangent bundle $T (\overline{V}^{n-1})$: 
$f_*: T (V^{n-1}) \rightarrow T (\overline{V}^{n-1})$ such that 
$f_*|_{V^{n-1}} = f$ and $f_*| T_x (V^{n-1})$ is 
a linear nondegenerate mapping.

We will prove now the theorem on conformal rigidity of 
hypersurfaces.

\noindent 
{\bf Theorem} {\em Let $n \geq 4$,  and  
 $V^{n-1}$ and $\overline{V}^{n-1}$ be two nonisotropic 
 hypersurfaces without umbilical points 
in a real conformal space $C^n$ or a real pseudoconformal space $C^n_q$. 
Suppose that there is a one-to-one correspondence 
$f: V^{n-1} \rightarrow \overline{V}^{n-1}$ between points of
 these hypersurfaces, 
and in the corresponding points of 
$V^{n-1}$ and $\overline{V}^{n-1}$ the following condition holds:
\begin{equation}\label{eq:20}
\overline{I} = f_* I,
\end{equation}
where $I$ and $\overline{I}$ are the conformal 
quadratic elements of 
the hypersurfaces $V^{n-1}$  and $\overline{V}^{n-1}$ 
defined above. 
 Then   the hypersurfaces $V^{n-1}$ and 
$\overline{V}^{n-1}$ are conformally equivalent, i.e. 
there exists a conformal transformation $\varphi$ of 
the space $C^n$ or $C^n_q$ 
such that  $\varphi (V^{n-1}) = \overline{V}^{n-1}$.}

{\sf Proof.} 
Relation (20) can be written in the form 
\begin{equation}\label{eq:21}
\frac{\overline{h}^2}{\overline{g}} = \frac{h^2}{g},  
\end{equation}
from which it follows that
$$
g \cdot \overline{h}^2 = \overline{g} \cdot h^2.
$$
Here and in what follows, for simplicity, we write $g$ and $h$ 
instead of $f_* g$ and $f_* h$. 
By the above lemma, the forms 
$ \overline{h}^2,  \overline{g}$  and $h^2, g$
do not have common factors. It follows that 
\begin{equation}\label{eq:22}
\overline{g} = \sigma^2 g, \;\; \overline{h} = \sigma h, 
\end{equation}
where $\sigma = \sigma (x) \neq 0$.

Further, let $x \in V^{n-1}$ and $\overline{x} 
\in \overline{V}^{n-1}$ be two corresponding points of the 
hypersurfaces $V^{n-1}$ and $\overline{V}^{n-1}$, and let 
$\varphi$ be a conformal transformation mapping $x = A_0$ into 
$\overline{x} = \overline{A}_0$ and the central tangent hypersphere $C_n$ 
into the central tangent  hypersphere $\overline{C}_n$.
Then the equations of $V^{n-1}$ and 
 $\overline{V}^{n-1}$  have the form:
\begin{equation}\label{eq:23}
 \omega_0^n = 0, \;\;\overline{\omega}_0^n =0.
\end{equation}
Moreover, the basis forms of $V^{n-1}$ and $\overline{V}^{n-1}$ are 
equal:
\begin{equation}\label{eq:24}
 \overline{\omega}_0^i = \omega_0^i.
\end{equation}
Since the first fundamental forms of 
the hypersurfaces $V^{n-1}$ and $\overline{V}^{n-1}$ 
have the form
$$
g = g_{ij} \omega^i \omega^j, \;\;\; \overline{g} = 
\overline{g}_{ij} \omega^i \omega^j,
$$
and their second fundamental forms have the form 
$$
h  = h_{ij} \omega^i \omega^j, \;\;\; 
\overline{h} = \overline{h}_{ij} \omega^i \omega^j,
$$
where $g^{ij} h_{ij} = 0$ and $\overline{g}^{ij} \overline{h}_{ij} 
= 0$ (see (13)), 
it follows that relations (22)  are equivalent to the relations
\begin{equation}\label{eq:.25}
 \overline{g}_{ij} = \sigma^2 g_{ij}
\end{equation}
and 
\begin{equation}\label{eq:26}
 \overline{h}_{ij} = \sigma h_{ij}.
\end{equation}

Since $\sigma \neq 0$, then by renormalizing the point $A_0$, 
this factor can be reduced to 1. In fact, setting
$$
 \overline{A}'_0 = \frac{1}{\sigma} \overline{A}_0,
$$
we find that
$$
d \overline{A}'_0 = d \Biggl(\frac{1}{\sigma}\Biggr) 
\overline{A}_0 + 
 \frac{1}{\sigma} (\overline{\omega}_0^0 \overline{A}_0 
+ \omega^i \overline{A}_i)
$$
and 
$$
(d \overline{A}'_0, d \overline{A}'_0) = \frac{1}{\sigma^2} 
 \overline{g}_{ij}  \omega^i \omega^j = g_{ij} \omega^i \omega^j.
$$
Thus,  we obtained 
\begin{equation}\label{eq:27}
 \overline{g}_{ij} = g_{ij}.
\end{equation}
After the above normalization,  
we have 
$$\overline{h}' = (d^2 \overline{A}_0', 
\overline{C}_n') = \frac{1}{\sigma} (d^2 \overline{A}_0, 
\overline{C}_n) = \frac{1}{\sigma} \overline{h} 
= h.
$$
It follows that  
\begin{equation}\label{eq:28}
 \overline{h}_{ij} = h_{ij}.
\end{equation}
Note that in (27) and (28) we wrote 
$\overline{g}_{ij}$ and  
$\overline{h}_{ij}$ instead of  
$\overline{g}_{ij}'$ and  
$\overline{g}_{ij}'$.

Taking the exterior derivatives of equations (24), we obtain
$$
[\overline{\omega}^i_j - \omega^i_j 
- \delta^i_j (\overline{\omega}^0_0 - \omega^0_0)] 
\wedge \omega^j = 0.
$$
Applying Cartan's lemma to these equations, we find that 
\begin{equation}\label{eq:29}
\overline{\omega}^i_j - \omega^i_j = 
 \delta^i_j (\overline{\omega}^0_0 - \omega^0_0) + T_{jk}^i 
\omega^k,
\end{equation}
where $ T_{jk}^i =  T_{kj}^i$. It is easy to prove that the 
quantities $T^i_{jk}$ form a $(1, 2)$-tensor, which is called the  
{\em deformation tensor} of the tangent bundle. 

Differentiating equations (27), we obtain
$$
g_{ik} (\overline{\omega}^k_j - \omega^k_j) 
+ g_{kj} (\overline{\omega}^k_i - \omega^k_i) = 0.
$$
Substituting for $\overline{\omega}^i_j - \omega^i_j$ 
 the values taken from (29), we find  that 
\begin{equation}\label{eq:30}
2 g_{ij} (\overline{\omega}^0_0 - \omega^0_0) 
+ (g_{ik} T^k_{jl} + g_{kj} T^k_{il}) \omega^l = 0.
\end{equation}
It follows that the 1-form $ \overline{\omega}^0_0 - \omega^0_0$ 
is expressed in terms of the basis forms $\omega^l$:
\begin{equation}\label{eq:31}
\overline{\omega}^0_0 - \omega^0_0 = s_l  \omega^l.
\end{equation}

Next, we  make the transformation 
$$
 \overline{A}'_i = \overline{A}_i + x_i \overline{A}_0
$$
in the pencil of normal hyperspheres. Since 
$ \overline{A}'_0 = \overline{A}_0$, we have 
$$
d \overline{A}'_0 = \overline{\omega}_0^0 \overline{A}_0 
+ \omega^i \overline{A}_i = \overline{\omega}_0^0 \overline{A}_0 
+ \omega^i ( \overline{A}'_i - x_i \overline{A}_0).
$$
It follows that
$$
{}'\overline{\omega}^0_0 =  \overline{\omega}_0^0 
- x_i \omega^i.
$$
By (31), from this we find that 
$$
{}'\overline{\omega}^0_0 =  \omega_0^0 + (s_i - x_i) \omega^i.
$$
We can see now that by setting $x_i = s_i$, we reduce relation 
(31) to the form
\begin{equation}\label{eq:32}
\overline{\omega}^0_0 = \omega^0_0.
\end{equation}

By (32), equations (30) take the form
$$
 g_{il} T^l_{jk} + g_{jl} T^l_{ik} = 0.
$$
By cycling these equations with respect to the indices $i, j$ and 
$k$ and subtracting the first equation from the sum of the last 
two equations, we obtain the conditions 
$$
T^k_{ij} = 0,
$$
by means of which equations (29) become 
\begin{equation}\label{eq:33}
\overline{\omega}^i_j = \omega^i_j.
\end{equation}

Taking the exterior derivatives of equations (32), we obtain 
the exterior quadratic equation
$$
(\overline{\omega}^0_i - \omega^0_i) \wedge \omega^i = 0,
$$
from which, by Cartan's lemma, it follows that
\begin{equation}\label{eq:34}
\overline{\omega}^0_i - \omega^0_i = t_{ij} \omega^j, \;\;
t_{ij} = t_{ji}.
\end{equation}

Taking the exterior derivatives of (33), we obtain 
$$
\omega^i \wedge  (\overline{\omega}^0_j - \omega^0_j) 
+ \overline{\omega}^i_n \wedge \overline{\omega}^n_j 
- \omega^i_n  \wedge \omega^n_j  
+ g^{ik} g_{jl} (\overline{\omega}^0_k - \omega^0_k) 
 \wedge \omega^l  = 0.
$$
By (28), the second and third terms on the left-hand side  
cancel out. Substituting for $\overline{\omega}^0_i - \omega^0_i$ 
in the remaining terms  the values taken from (34), and 
using the fact that the forms $\omega^l$ are linearly 
independent, we find that 
$$
- t_{jk} \delta^i_l +  t_{jl} \delta^i_k + g^{im} 
(t_{mk} g_{jl} -  t_{ml} g_{jk}) = 0.
$$
Contracting this relation with respect to the indices $i$ and 
$k$, we arrive at the equation
\begin{equation}\label{eq:35}
  (n - 3) t_{jl} = - t g_{jl},
\end{equation}
where $t = g^{im} t_{im}$. Since $n \geq 4$, 
by contracting the latter equation with 
the tensor $g^{jl}$, we find that $(2n - 4) t  = 0$.  
It follows  that $t = 0$, and consequently $t_{jl} = 0$.

As a result, equation (34) takes the form
\begin{equation}\label{eq:36}
\overline{\omega}^0_i = \omega^0_i.
\end{equation}

 From equations (14) it follows that 
\begin{equation}\label{eq:37}
\overline{\omega}^n_i = \omega^n_i.
\end{equation}
Taking the exterior derivatives of equations (37), we 
obtain 
\begin{equation}\label{eq:38}
\omega_i^{n+1} \wedge (\overline{\omega}^0_n - \omega^0_n) = 0.
\end{equation}
By (7), even for $n \geq 3$, the forms $\omega_i^{n+1}$ are linearly 
independent, and as a result, we find that 

\begin{equation}\label{eq:39}
\overline{\omega}^0_n = \omega^0_n.
\end{equation}
Exterior differentiation of equations (36) and (39) leads to 
the identity. Thus the system of equations 
 (23), (24), (27), (32), (33), (36) (37) and (39)
is completely integrable. 

Moreover, equations  
 (23), (24), (27), (32), (33), (36) (37) and (39)
show that all components of an 
infinitesimal displacement of second order moving 
frames  associated with the hypersurfaces $V^{n-1}$ and 
$\overline{V}^{n-1}$ coincide. Thus,   by the equivalence 
theorem  of \'{E}. Cartan (see [C 08] or [Ga 89]),  
the hypersurface $\overline{V}^{n-1}$ can be obtained from the 
hypersurface $V^{n-1}$ by means of a conformal transformation. 
Therefore,  the hypersurfaces $V^{n-1}$ and $\overline{V}^{n-1}$ 
are conformally equivalent. \rule{3mm}{3mm}

As we can see from  equation (35), 
the proof of our main theorem  
fails if $n = 3$. To prove the rigidity theorem 
for $n = 3$,  it is necessary to add 
certain additional conditions to condition (20) which are 
connected with a third order differential neighborhood. 
Such a theorem for $n = 3$ was proved in  [SSu 80].

{\sf             
Department of Mathematics,            
Ben-Gurion University of the Negev,  
P.O. Box 653,  Beer-Sheva 84105, Israel       
}

{\em E-mail address}: akivis@black.bgu.ac.il

\vspace*{3mm}

{\sf 
 Department of Mathematics,   New Jersey Institute of Technology,  
 University Heights, Newark, NJ 07102
}

{\em E-mail address}: vlgold@numerics.njit.edu


\begin{thebibliography}{BCGGG 91}

\bibitem[AG 93]{AG:[AG 93]} M. A. Akivis and V. V. Goldberg,  
{\em Projective differential geometry of submanifolds}. 
North-Holland, Amsterdam-New York-Tokyo, 1993, xi+362 pp. 


\bibitem[CD 87]{CD:[CD 87]} M. do Carmo and M. Dajczer,  
{\em Conformal rigidity}. Amer. J. Math. {\bf 109} (1987), no. 5,
 963--985. 

\bibitem[C 08]{C:[Ca 08]} \'{E}. Cartan,  
{\em Les sous-groupes des groupes continus de transformations.} 
Ann. Sci. \'{E}cole Norm. (3) {\bf 25} (1908), 57--194. 

\bibitem[C 17]{C:[C 17]} \'{E}. Cartan, 
{\em La d\'{e}formation des hypersurfaces dans l'espace conforme
 r\'{e}el \`{a} $n \geq 5$ dimensions}.  Bull. Soc. Math. France
 {\bf 45} (1917), 57--121.


\bibitem[Ca 20]{C:[C 20]}
\'{E}. Cartan,  
 {\em Sur la d\'{e}formation projective des surfaces}.  Ann. Sci.
 \'{E}cole Norm. Sup. {\bf 37} (1920), 259--356.  


\bibitem[F 16]{F:[F 16]} G.  Fubini,  
{\em Applicabilit\'{a} proiettiva di due superficie}. Rend. Circ. 
Mat Palermo {\bf 41} (1916), 135--162.


 \bibitem[F 18]{F:[F 18]}  G.  Fubini,  
{\em Stud\^{\i} relativi all' elemento lineare proiettivo di una 
ipersuperficie}. Atti Accad. Naz. Lincei Rend. Cl. Sci. Fis. 
Mat. Natur. (5) {\bf 27} (1918), 99--106.


 \bibitem[F\v{C} 26]{FC:[FC 26]} G.  Fubini, and E.  \v{C}ech,  
{\em Geometria proiettiva differenziale}. Zanichelli, Bologna, 
vol. 1, 1926, 394 pp., vol. 2, 1927,  400 pp. 


\bibitem[Gar 89]{Gar:[Gar 89]} R. Gardner,   
{\em The method of equivalence and its applications}. 
CBMS-NSF Regional conference Series in Applied 
Mathematics, {\bf 58}. SIAM, Philadelphia, PA, 1989, vii+127 pp. 


\bibitem[JM 92]{JM:[JM 92]}   G. R.  Jensen and E. Musso, 
{\em Rigidity of hypersurfaces in complex projective space}, 
Ann. Sci. \'{E}cole Norm. Sup. (4) {\bf 27} (1994), 227--248.


\bibitem[S 62]{S:[S 62]}   R. Sacksteder,  
{\em The rigidity of hypersurfaces}, 
J. Math. Mekh.  {\bf 11} (1962), 929--940.

\bibitem[SSu 80]{SSu:[SSu 80]} C. Schiemankgk and R. Sulanke,  
{\em Submanifolds of the M\"{o}bius space}, Math. Nachr. 
 {\bf 96} (1980), 165--183. 

\bibitem[Su 82]{Su:[Su 82]}  R. Sulanke, 
{\em  Submanifolds of the M\"{o}bius space. III The analogue of
 O. Bonnet's theorem for hypersurfaces}. Tensor (N. S.) {\bf 38}  
(1982), 311--317.

\end{thebibliography}
 \end{document}